\documentclass[12pt]{amsart}
\usepackage{amsmath, amssymb, euscript, hyperref}

\theoremstyle{definition}

\begin{document}

\address{Azer Akhmedov, Department of Mathematics,
North Dakota State University,
Fargo, ND, 58108, USA}
\email{azer.akhmedov@ndsu.edu}

 \begin{center} {\bf \Large A weak Zassenhaus lemma for discrete subgroups of Diff(I)} \end{center}
  
  \bigskip

 \begin{center} {\bf Azer Akhmedov} \end{center}
 
 \vspace{1cm}
 
 ABSTRACT: {\Small We prove a weaker version of the Zassenhaus Lemma for subgroups of Diff(I). We also show that a group with commutator subgroup containing a non-abelian free subsemigroup does not admit a $C_0$-discrete faithful representation in Diff(I).}

\vspace{1cm}

  In this paper, we continue our study of discrete subgroups of $\mathrm{Diff} _{+}(I)$ - the group of orientation preserving diffeomorphisms of the closed interval $I = [0,1]$. Following recent trends, we try to view the group $\mathrm{Diff} _{+}(I)$ as an analog of a Lie group, and we study still basic questions about discrete subgroups of it. This paper can be viewed as a continuation of {\bf [A]} although the proofs of the results of this paper are independent of {\bf [A]}. 
  
  \medskip   
  
   Throughout the paper, the letter $G$ will denote the group $\mathrm{Diff} _{+}(I)$. On $G$, we assume the metric induced by the standard norm of the Banach space $C^1[0,1]$. We will denote this metric by $d_1$. Sometimes, we also will consider the metric on $G$ that comes from the standard sup norm $||f||_0 = \mathrm{sup}\ _{x\in [0,1]}|f(x)|$ of $C[0,1]$ which we will denote by $d_0$. However, unless specified, the metric in all the groups $\mathrm{Diff} _{+}^{r}(I), r\in \mathbb{R}, r \geq 1$ will be assumed to be $d_1$.
 
 \medskip
 
  The central theme of the paper is the Zassenhaus Lemma. This lemma states that in a connected Lie group $H$ there exists an open non-empty neighborhood $U$ of the identity such that any discrete subgroup generated by elements from $U$ is nilpotent (see {\bf [R]}). For example, if $H$ is a simple Lie group (such as $SL_2(\mathbb{R})$), and $\Gamma \leq H$ is a lattice, then $\Gamma $ cannot be generated by elements too close to the identity.  
     
     \bigskip
     
  In this paper we prove weak versions of the Zassenhaus Lemma for the group $G = \mathrm{Diff} _{+}(I)$. Our study leads us to showing that finitely generated groups with exponential growth which satisfy a very mild condition do not admit faithful $C_0$-discrete representation in $G$:  
     
     \medskip
     
     {\bf Theorem A.} Let $\Gamma $ be a subgroup of $G$, and $f, g\in [\Gamma ,\Gamma ]$ such that $f$ and $g$ generate a non-abelian free subsemigroup. Then $\Gamma $ is not $C_0$-discrete. 
     
     \medskip
     
     We also study the Zassenhaus Lemma for the relatives of $G$ such as $\mathrm{Diff} _{+}^{1+c}(I), c\in \mathbb{R}, c > 0$ - the group of orientation preserving diffeomorphisms of regularity $1+c$. In the case of $\mathrm{Diff} _{+}^{1+c}[0,1]$, combining Theorem A with the results of {|bf [N2]}, we show that $C_0$-discrete subgroups are more rare. 
     
     \medskip
     
     {\bf Theorem B.} Let $\Gamma $ be a $C_0$-discrete subgroup of $\mathrm{Diff} _{+}^{1+c}[0,1]$. Then $\Gamma $ is solvable with solvability degree at most $k(c)$.
     
     \medskip
     
     Theorem B can be strengthened if the regularity is increased further; combining Theorem A with the results of Navas {\bf [N2]}, Plante-Thurston {\bf [PT]}, and Szekeres {\bf [S]} we obtain the following
     
     \medskip
     
     {\bf Theorem C.} If $\Gamma $ is $C_0$-discrete subgroup of $\mathrm{Diff} _{+}^{2}[0,1]$ then $\Gamma $ is metaabelian.
    
     \bigskip
     
     It follows from the results of {\bf [A]}, as remarked there, that the Zassenhaus Lemma does not hold either for $\mathrm{Diff} _{+}(I)$ or for $\mathrm{Homeo} _{+}(I)$, in metrics $d_1$ and $d_0$ respectively. 
     
     \medskip
     
   In the increased regularity the lemma still fails: given an arbitrary open neighborhood $U$ of the identity diffeomorphism in $G$, it is easy to find two $C^{\infty }$ ``bump functions" in $U$ which generate a discrete group isomorphic to $\mathbb{Z}\wr \mathbb{Z}$; thus the lemma fails for $\mathrm{Diff} _{+}^{\infty }(I)$. 
   
    \medskip
    
   Because of the failure of the lemma, it is natural to consider strongly discrete subgroups which we have defined in {\bf [A]}. Indeed, for strongly discrete subgroups, we are able to obtain positive results which are natural substitutes for the Zassenhaus Lemma.

    \medskip
    
  Let us recall the definition of strongly discrete subgroup from {\bf [A]}: 
  
 \medskip
 
   {\bf Definition 1.} Let $\Gamma $ be a subgroup of $\mathrm{Diff} _{+}(I)$. $\Gamma $ is called {\em strongly discrete} if there exists $C>0$ and $x_0\in (0,1)$ such that $|g'(x_0)-1| > C$ for all $g\in \Gamma \backslash \{1\}$.  Similarly, we say $\Gamma $ is $C_0$-{\em strongly discrete} if $|g(x_0)-x_0| > C$ for all $g\in \Gamma \backslash \{1\}$.

 \bigskip
  
  Let us note that a strongly discrete subgroup of $G$ is discrete; and $C_0$-strongly discrete subgroup of $G$ is $C_0$-discrete.
  
  \medskip
  
  For convenience of the reader, let us recall several basic notions on the growth of groups: if $\Gamma $ is a finitely generated group, and $S$ a finite generating set, we will define $\omega (\Gamma , S) = \lim _{n\rightarrow \infty } \sqrt[n]{|B_n(1; S, \Gamma )|}$ where $B_n(1; S, \Gamma )$ denotes the ball of radius $n$ around the identity element. (Often we will denote this ball simply by $B_n(1)$). We will also write $\omega (\Gamma ) = \mathrm{inf} _{|S| < \infty , \langle S\rangle =\Gamma}\omega (\Gamma ,S)$ where the infimum is taken over all finite generating sets $S$ of $\Gamma $. If $\omega (\Gamma ) > 1$ then one says that $\Gamma $ has uniform exponential growth. 
  
  \medskip
  
  Now we are ready to state weak versions of the Zassenhaus Lemma for the group $G$. First, we state a theorem about $C_0$-strongly discrete subgroups.
 
 \medskip
 
 {\bf Theorem 1.} Let $\omega > 1$. Then there exists an open non-empty neighborhood $U$ of the identity $1\in \mathrm{Diff} _{+}^{1}[0,1]$ such that if $\Gamma $ is a finitely generated $C_0$-strongly discrete subgroup of $\mathrm{Diff} _{+}^{1}[0,1]$ with $\omega (\Gamma ) \geq \omega $, then $\Gamma $ cannot be generated by elements from $U$.

 \bigskip
 
  By increasing the regularity, we can prove a similar version for strongly discrete subgroups
  
  \medskip
     
  {\bf Theorem 2.} Let $\omega > 1$. Then there exists an open non-empty neighborhood $U$ of the identity $1\in \mathrm{Diff} _{+}^{1}[0,1]$ such that if $\Gamma $ is a finitely generated  strongly discrete subgroup of $\mathrm{Diff} _{+}^{2}[0,1]$ with $\omega (\Gamma ) \geq \omega $, then $\Gamma $ cannot be generated by elements from $U$. 
  
    \bigskip
  
  {\bf Remark 1.} In regard to the Zassenhaus Lemma, it is interesting to ask a reverse question, i.e. given an arbitrary open neighborhood $U$ of the identity in $G$, is it true that any finitely generated torsion free nilpotent group $\Gamma $ admits a faithful discrete representation in $G$ generated by elements from $U$? In {\bf [FF]}, it is proved that any such $\Gamma $ does admit a faithful representation into $G$ generated by diffeomorphisms from $U$. Also, it is proved in {\bf [N1]} that any finitely generated nilpotent subgroup of $G$ indeed can be conjugated to a subgroup generated by elements from $U$.     
  
  \medskip
  
  {\bf Remark 2.} Because of the assumptions about uniform exponential growth in Theorem 1 and Theorem 2, it is natural to ask whether or not every finitely generated subgroup of $G$ of exponential growth has uniformly exponential growth.
This question has already been raised in {\bf [N2]}.

  \bigskip
  
  {\em Acknowledgement.} We are thankful to Andr\'es Navas for his remarks, and for bringing {\bf [N1]} to our attention. We also would like to thank to two anonymous referees for many helpful remarks and suggestions.

  \vspace{1cm}
  
 \begin{center} {\Large Proofs of Theorem 1 and Theorem 2} \end{center}
   
  \medskip
  
{\bf Proof of Theorem 1.} We can choose $\lambda > 1$ such that $\lambda < \omega (\Gamma )$. Then the cardinality of the sphere of radius $n$ of $\Gamma $ with respect to any fixed finite generating set is bigger than the exponential function $\lambda ^n$, for infinitely many $n$.

 \medskip

 Then let $\epsilon > 0$ such that $(1-10\epsilon)\lambda > 1$. We let $U$ be the $\epsilon $-neighborhood of the identity in $G$ with respect to $d_1$ metric (we always assume $d_1$ metric in $G$ unless otherwise stated).
  
  \medskip
  
  Let $\Gamma $ be generated by finitely many non-trivial diffeomorphisms $f_1, f_2, \ldots ,f_s\in U$. We fix this generating set and denote it by $S$, i.e. $S = \{f_1, f_1^{-1}, \ldots , f_s, f_s^{-1}\}$.

 \medskip
 
  We want to prove that $\Gamma $ is not $C_0$-strongly discrete. Assuming the opposite, let  $x_0\in (0,1)$ such that for some $C > 0$, \ $|g(x_0)-x_0| > C$ for all $g\in \Gamma \backslash \{1\}$. 
  
  \medskip
  
  Let $B_n(1)$ be the ball of radius $n$ around the identity in the Cayley graph of $\Gamma $ with respect to $S$. Then  $\mathrm{Card}(B_n(1)\backslash B_{n-1}(1)) > \lambda ^n$ for infinitely many $n\in \mathbb{N}$. Let $A$ denotes the set of all such $n$.
  
\medskip

   Let $\Delta $ be a closed subinterval of $(0,1)$ of length less than $C$ such that $x_0$ is the left end of $\Delta $.

 \medskip
   
 We denote the right-invariant Cayley metric of $\Gamma $ with respect to $S$ by $|.|$. For all $g\in \Gamma $, let $\Delta _g = g(\Delta )$. Thus we have a collection $\{\Delta _g\}_{g\in G}$ of closed subintervals of $(0,1)$.
  
 \medskip

 Notice that if $g = sw, s\in S$ then by Mean Value Theorem, $|\Delta _{sw}| > (1-10\epsilon )|s(\Delta _{w})|$. Then, necessarily, for all $n\in A$, we have $\sum _{|g| = n}|\Delta _g| > (1-10\epsilon )^n\lambda ^n|\Delta | \rightarrow \infty $ as $n\rightarrow \infty $.
 
 \medskip
 
 Then $\exists g_1, g_2\in \Gamma , g_1\neq g_2$ such that $g_2(x_0) \in \Delta _{g_1}$. Then $g_1^{-1}g_2(x_0)\in \Delta $. Since $|\Delta | < C $ we obtain a contradiction. \ $\square $
 
 \bigskip
  
 Now we prove a better result by assuming higher regularity for the representation
 
 \medskip
 
 {\bf Proof of Theorem 2.} Let $\lambda , \lambda _1, \lambda _2 $ be  constants such that $1 < \lambda < \lambda _1 < \lambda _2 < \omega (\Gamma )$. Then the cardinality of the sphere of radius $n$ of $\Gamma $ with respect to any fixed finite generating set is bigger than the exponential function $\lambda _2^n$, for infinitely many $n$. 
 
 \medskip
 
 We choose $\epsilon > 0, \eta > 0$ to be such that $1 < \eta < \frac{\lambda }{1+\epsilon }$ and $\frac {1+\epsilon }{1-\epsilon } < \frac{\lambda _1}{\lambda }$. Let $U$ be the ball of radius $\epsilon $ around the identity diffeomorphism.
 
 \medskip
 
  We again assume $\Gamma $ is generated by finitely many non-trivial diffeomorphisms $f_1, f_2, \ldots ,f_s\in U$, and we fix the generating set $S = \{f_1, f_1^{-1}, \ldots , f_s, f_s^{-1}\}$. Let $B_n(1)$ be the ball of radius $n$ around the identity in the Cayley graph of $\Gamma $ with respect to $S$. Then we have $\mathrm{Card}(B_n(1)\backslash B_{n-1}(1)) > \lambda _2^n$ for infinitely many $n\in \mathbb{N}$. Let $A$ denotes the set of all such $n$.
 
 \medskip
 
  We need to show that $\Gamma $ is not strongly discrete. Assuming the opposite, let $x_0\in (0,1)$ such that for some $C > 0$, \ $|g'(x_0)-1| > C$ for all $g\in \Gamma \backslash \{1\}$. 
 
 \medskip
 
 Let $C_1$ be a positive number such that  $$1-C < (1-C_1)^2 \ \mathrm{and}  \ 1+C > (1+C_1)^2$$ Let also $N_1\in \mathbb{N}$ such that for any $n \geq N_1$, we have $$1-C_1 < (1-\frac{1}{\eta ^n})^{n} < (1+\frac{1}{\eta ^n})^{n} < 1 + C_1  $$
 
 \medskip
 
 Notice that for all $n\in A, g\in B_n(1)\backslash B_{n-1}(1)$, and $x\in [0,1]$, we have $(1-\epsilon )^n < g'(x) < (1+\epsilon )^n$. Since $\frac {1+\epsilon }{1-\epsilon } < \frac{\lambda _1}{\lambda }$, there exists $n\in A$ and $g_1, g_2\in \Gamma $ such that $$n > N_1, \ g_1\neq g_2, \ |g_1| = |g_2| = n$$ but $$|g_1(x_0)-g_2(x_0)| \leq \frac{1}{\lambda ^n} \ (\star _1), \  \mathrm{and} \ 1 - C_1  < \frac {g_1'(x_0)}{g_2'(x_0)} < 1 + C_1 \ (\star _2)  $$
 
 \medskip
 
  Indeed, by the pigeonhole principle, for all $n\in A$, there exists $j\in \{0, 1, \ldots , [\lambda ^n]\}$ such that $$\mathrm{Card}\{g\in B_n(1)\backslash B_{n-1}(1) \ | \ g(x_0)\in [\frac{j}{\lambda ^n}, \frac{j+1}{\lambda ^n})\}\geq \frac{\lambda _2^n}{{\lambda ^n}+1}$$ 
 
 \medskip
 
 For all $n\in A, j\in \{0, 1, \ldots , [\lambda ^n]\}$, let $$D(n,j) = \{g\in B_n(1)\backslash B_{n-1}(1) \ | \ g(x_0)\in [\frac{j}{\lambda ^n}, \frac{j+1}{\lambda ^n})\}$$ Then, for sufficiently big $n\in A$, there exists $j\in \{0, 1, \ldots , [\lambda ^n]\}$ such that $\mathrm{Card}(D(n,j)) \geq \frac{\lambda _1^n}{\lambda ^n} \ (\star _3)$. For all $n\in A$, let $$J(n) = \{j\in \{0, 1, \ldots , [\lambda ^n]\} \ | \ \mathrm{Card}(D(n,j)) \geq \frac{\lambda _1^n}{\lambda ^n}\}$$  
 
 \medskip
 
 Recall also that for all $g\in D(n,j)$, we have $$(1-\epsilon )^n < g'(x_0) < (1+\epsilon )^n$$ Then, since  $\frac {1+\epsilon }{1-\epsilon } < \frac{\lambda _1}{\lambda }$, for sufficiently big $n\in A$ and $j\in J(n)$, applying the pigeonhole principle to the set $D(n,j)$, we obtain that (besides the inequality $(\star _3)$) there exist distinct  $g_1, g_2\in D(n,j)$ such that the inequality $1 - C_1  < \frac {g_1'(x_0)}{g_2'(x_0)} < 1 + C_1$ holds. On the other hand, by definition of $D(n,j)$, we have $|g_1(x_0)-g_2(x_0)| \leq \frac{1}{\lambda ^n}$; thus we established the desired inequalities $(\star _1)$ and $(\star _2)$.
 
 \medskip
 
 Let now $y_0 = g_1(x_0), z_0 = g_2(x_0), W = g_1^{-1}, V = g_1^{-1}g_2$, and let $W = h_nh_{n-1}\ldots h_1$ where $W$ is a reduced word in the alphabet $S$ of length $n$ and $h_i\in S, 1\leq i\leq n$. 
 
 \medskip
 
 Let also $W_k$ be the suffix of $W$ of length $k, \ y_k = W_k(y_0), z_k = W_k(z_0), \ 1\leq k\leq n$.
 
 \medskip
 
 Furthermore, let $\mathrm{max} _{1\leq i\leq s} \mathrm{sup}  _{0\leq y\neq z\leq 1}\frac{|f_i'(y)-f_i'(z)|}{|y-z|} = M$, and $L = 1 + \epsilon $. 
    
 \medskip
 
 Then we have $$|y_k-z_k|\leq \frac{L^k}{\lambda ^n}, \ |h_{k+1}'(y_k)-h_{k+1}'(z_k)|\leq \frac{ML^k}{\lambda ^n}, 0\leq k\leq n-1$$
 
 \medskip
 
    Then $1 - \frac {ML^{k+1}}{\lambda ^n} \leq \frac {h_{k+1}'(y_k)}{h_{k+1}'(z_k)} \leq 1 + \frac {ML^{k+1}}{\lambda ^n}$ for all $0\leq k\leq n-1$. From here we obtain that
   $$\Pi _{k=0}^{n-1}(1 - \frac {ML^{k+1}}{\lambda ^n})\leq \Pi _{k=0}^{n}\frac{h_{k+1}'(y_k)}{h_{k+1}'(z_k)}\leq \Pi _{k=0}^{n-1}(1 + \frac {ML^{k+1}}{\lambda ^n})$$
   
   \medskip
   
   Then, for sufficiently big $n$ in $A$  $$(1 - \frac {1}{\eta ^n})^n = \Pi _{k=0}^{n-1}(1 - \frac {1}{\eta ^n})\leq \Pi _{k=0}^{n-1}\frac{h_{k+1}'(y_k)}{h_{k+1}'(z_k)}\leq \Pi _{k=0}^{n-1}(1 + \frac {1}{\eta ^n}) = (1 + \frac {1}{\eta ^n})^n$$
   
   \medskip
   
   Since, by the chain rule, $\Pi _{k=0}^{n-1}\frac{h_{k+1}'(y_k)}{h_{k+1}'(z_k)} = \frac {(g_1^{-1})'(y_0)}{(g_1^{-1})'(z_0)}$, we obtain that $$1 - C_1 < \frac {(g_1^{-1})'(y_0)}{(g_1^{-1})'(z_0)} < 1 + C_1 $$

   \medskip
   
   Then 
   
   $$V'(x_0) = (g_1^{-1})'(g_2(x_0))g_2'(x_0) = \frac {(g_1^{-1})'(g_2(x_0))g_2'(x_0)}{(g_1^{-1})'(g_1(x_0))g_1'(x_0)}(g_1^{-1})'(g_1(x_0))g_1'(x_0) = $$ \ $$\frac {(g_1^{-1})'(g_2(x_0))g_2'(x_0)}{(g_1^{-1})'(g_1(x_0))g_1'(x_0)} = \frac {(g_1^{-1})'(g_2(x_0))}{(g_1^{-1})'(g_1(x_0))} \frac {g_2'(x_0)}{g_1'(x_0)} \in ((1-C_1)^2, (1 + C_1)^2) \subset (1-C, 1+C)$$

   \medskip
     
   Thus we proved that $1-C < V'(x_0) < 1 + C$ which contradicts our assumption. $\square $
   
   \bigskip
   
  {\bf Remark 3.} The same proof, with slight changes, works for representations of $C^{1+c}$-regularity for any real $c>0$.
  
  \vspace{1cm}
  
 \begin{center} {\Large Proofs of Theorems A, B, C} \end{center}
  
  \bigskip
  
  In the proofs of Theorem 1 and of Theorem 2, we consider the orbit of the point $x_0$ under the action of $\Gamma $. By using exponential growth, we find two distinct elements $g_1, g_2$ such that $g_1(x_0)$ and $g_2(x_0)$ are very close. Then we ``pull back" $g_2(x_0)$ by $g_1^{-1}$, i.e. we consider the point $g_1^{-1}g_2(x_0)$ and show that this point is sufficiently close to $x_0$. It is at this stage that we heavily use the condition that $\Gamma $ is generated by elements from the small neighborhood of $1\in G$, i.e. derivatives of the generators are uniformly close to 1. However, if $\Gamma $ is an arbitrary subgroup of the commutator group $[G,G]$, not necessarily generated by elements close to the identity element, then for any $x_0\in (0,1), f\in \Gamma $ and for any $\epsilon > 0$, there exists $W\in \Gamma $ such that $|f'(W(x_0))-1|< \epsilon $; we simply need to find $W$ such that $W(x_0)$ is sufficiently close to 1 (or to 0). This fact provides a new idea of taking $x_0$ close to 1, then considering the part of the orbit which lies in a small neighborhood of 1, then using exponential growth to find points close to each other in that neighborhood, and then perform the ``pull back".
  
  \medskip
  
  The following proposition is a special case of Theorem A, and answers Question 2 from {\bf [A]}. For simplicity, we give a separate proof of it.
  
  \medskip
  
  {\bf Proposition 1.} $\mathbb{F}_2$ does not admit a faithful  $C_0$-discrete representation in $G$.
  
  \medskip
  
  {\bf Proof.} Since the commutator subgroup of $\mathbb{F}_2$ contains an isomorphic copy of $\mathbb{F}_2$, it is sufficient to prove that $\mathbb{F}_2$ does not admit a faithful $C_0$-discrete representation in $G^{(1)} = [G,G]$. 
 
 \medskip
  
  Let $\Gamma $ be a subgroup of $G^{(1)}$ isomorphic to $\mathbb{F}_2$ generated by diffeomorphisms $f$ and $g$. Without loss of generality we may assume that $\Gamma $ has no fixed point on $(0,1)$. Let also $\epsilon > 0$ and $M = \mathrm{max} _{0\leq x\leq 1}(|f'(x)|+|g'(x)|)$. 
 
 \medskip
  
 We choose $N\in \mathbb{N}, \delta > 0$ and $\theta _N$ such that $1/N <  \epsilon , 1 < \theta _N < \sqrt[2N]{2}$, and for all $x\in [1-\delta , 1]$, the inequality $\frac{1}{\theta _N} < \phi '(x) < \theta _N$ holds where $\phi \in \{f, g, f^{-1}, g^{-1}\}$.
 
 \medskip
 
 Let $W = W(f,g)$ be an element of $\Gamma $ such that $W(1/N)\in [1-\delta , 1]$, $m$ be the length of the reduced word $W$. Let also $x_i = i/N, 0\leq i\leq N$.     
  
  \medskip
  
 For every $n\in \mathbb{N}$, let $$S_n = \{H\in B_n(1) \ | \ u(W(x_1))\geq W(x_1) \ \mathrm{for \ all \ suffixes} \ u \ \mathrm{of} \ H\}$$ (Here we view $H$ as a reduced word in the alphabet $\{f,g,f^{-1},g^{-1}\}$). Then $|S_n| \geq 2^n$. 
 
 \medskip
 
 Then (assuming $N\geq 3$) we can choose and fix a sufficiently big $n$ such that the following two conditions are satisfied:
 
 \medskip
 
 (i) there exist $g_1, g_2\in S_n$ such that $g_1\neq g_2$, and $$|g_1W(x_i)-g_2W(x_i)| < \frac{1}{\sqrt[2N]{2}^n}, 1\leq i\leq N-1.$$
 
 \medskip
 
 (ii) $M^m(\theta _N)^n\frac{1}{\sqrt[2N]{2}^n} < \epsilon $.
 
 \medskip
 
  Indeed, let $(c_0, c_1, \ldots , c_{N-1}, c_N)$ be a sequence of real numbers such that $\sqrt[2N]{2} = c_N < c_{N-1} < \ldots < c_1 < c_0 = 2$ and $c_i > \sqrt[2N]{2}c_{i+1}$, for all $i\in \{0, 1, \ldots , N-1\}$. Then, by the pigeonhole principle, for sufficiently big $n$, there exists a subset $S_n(1)\subseteq S_n$ such that $|S_n(1)| \geq c_1^n$ and $|g_1W(x_1)-g_2W(x_1)| < \frac{1}{\sqrt[2N]{2}^n}, \forall g_1, g_2\in S_n(1)$. 
  
  \medskip
  
  Suppose now $1\leq k\leq N-2$, and $S_n\supseteq S_n(1)\supseteq \ldots \supseteq S_n(k)$ such that for all $j\in \{1, \ldots , k\}$, $|S_n(j)|\geq c_j^n$ and for all $g_1, g_2\in S_n(j)$ we have $$|g_1W(x_i)-g_2W(x_i)| < \frac{1}{\sqrt[2N]{2}^n}, 1\leq i\leq j.$$
  
  Then by applying the pigeonhole principle to the set $S_n(k)$ for sufficiently big $n$, we obtain $S_n(k+1)\subseteq S_n(k)$ such that $|S_n(k+1)|\geq c_{k+1}^n$, and for all $g_1, g_2\in S_n(k+1)$ we have $$|g_1W(x_i)-g_2W(x_i)| < \frac{1}{\sqrt[2N]{2}^n}, 1\leq i\leq k+1.$$ 
  
 Then, for $k = N-2$, we obtain the desired inequality (condition (i)). 
 
 \bigskip
 
 Now, let $$h_1 = g_1W, h_2 = g_2W, y_i = W(x_i), z_i' = g_1(y_i), z_i'' = g_2(y_i), 1\leq i\leq N.$$ 
 
 \medskip
 
 Without loss of generality, we may also assume that $g_2(y_1) \geq g_1(y_1)$.
 
 \medskip
 
 Then for all $i\in \{1, \ldots , N-1\}$, we have
 
 \medskip
 
  $$|h_1^{-1}h_2(x_i) - x_i| = |(g_1W)^{-1}(g_2W)(x_i) - x_i| = $$ \ $$ |(g_1W)^{-1}(g_2W)(x_i) - (g_1W)^{-1}(g_1W)(x_i)| = |W^{-1}g_1^{-1}g_2(y_i) - W^{-1}g_1^{-1}g_1(y_i)| $$ \ $$=|
 W^{-1}g_1^{-1}(z_i'') - W^{-1}g_1^{-1}(z_i')|$$
 
 \medskip
 
  Let $u$ be a prefix of the reduced word $g_1$, and $g_1 = uv$ (so a reduced word $v$ is a suffix of $g_1$). Then, since $g_1, g_2\in S_n$, we have $$u^{-1}(z_i') = v(y_i)\geq v(y_1) \geq y_1$$ and $$u^{-1}(z_i'') = u^{-1}(g_2(y_i)) \geq u^{-1}(g_2(y_1)) \geq u^{-1}(g_1(y_1)) \geq v(y_1) \geq y_1$$ 
  
  \medskip
  
  Then by the Mean Value Theorem, we have $$|h_1^{-1}h_2(x_i) - x_i|\leq M^m(\theta _N)^n|z_1'-z_1''| < M^m(\theta _N)^n\frac{1}{\sqrt[2N]{2}^n}$$ Then, by condition (ii), we obtain $|h_1^{-1}h_2(x_i) - x_i| < \epsilon $. Then we have $|h_1^{-1}h_2(x) - x| < 2\epsilon $ for all $x\in [0,1]$. Indeed, let $x\in [x_i,x_{i+1}]$. Then $$|h_1^{-1}h_2(x) - x| \leq \mathrm{max}\{|h_1^{-1}h_2(x_i) - x|, |h_1^{-1}h_2(x_{i+1}) - x|\}$$
  
  \medskip
  
  But $|h_1^{-1}h_2(x_i) - x| \leq |h_1^{-1}h_2(x_i) - x_i| + |x_i-x| < 2\epsilon $, and similarly, $|h_1^{-1}h_2(x_{i+1}) - x| \leq |h_1^{-1}h_2(x_{i+1}) - x_{i+1}| + |x_{i+1}-x| < 2\epsilon $.
  
  \medskip
   
  Since $\epsilon $ is arbitrary, we obtain that $\Gamma $ is not $C_0$-discrete. $\square $       
  
  \vspace{1cm}

  By examining the proof of Proposition 1, we will now prove Theorem A thus obtaining a much stronger result. The inequality $|S_n| \geq 2^n$ is a crucial fact in the proof of Proposition 1; we need the cardinality of $S_n$ grow exponentially. If $\Gamma $ is an arbitrary finitely generated group with exponential growth, this exponential growth of $S_n$ is not automatically guaranteed. But we can replace $S_n$ by another subset $\mathbb{S}_n $ which still does the job of $S_n$ and which grows exponentially, if we assume a mild condition on $\Gamma $. 
    
  \medskip
  
  First we need the following easy
  
  \medskip
  
  {\bf Lemma 1.} Let $\alpha , \beta \in G , z_0\in (0,1)$ such that $z_0 \leq \alpha (z_0)\leq \beta \alpha (z_0)$. Then $U\beta \alpha (z_0)\geq z_0$ where $U = U(\alpha ,\beta )$ is any positive word in letters $\alpha , \beta $. $\square $  
    
  \bigskip
    
  Now we are ready to prove Theorem A.

  \medskip
  
  {\bf Proof.} Without loss of generality, we may assume that $\Gamma $ has no fixed point on $(0,1)$. Let again $\epsilon > 0, N\in \mathbb{N}, \delta > 0, \theta _N > 0,  M = 2\mathrm{sup} _{0\leq x\leq 1}(|f'(x)|+|g'(x)|)$ such that $1/N < \epsilon , 1 < \theta _N < \sqrt[2N]{2}$, and for all $x\in [1-\delta , 1]$, the inequality $\frac{1}{\theta _N} < \phi '(x) < \theta _N$ holds where $\phi \in \{f, g, f^{-1}, g^{-1}\}$.

 \medskip
 
 Let $W = W(f,g)$ be an element of $\Gamma $ such that $$(\{f^iW(1/N) \ | \ -2\leq i\leq 2\}\cup \{g^iW(1/N) \ | \ -2\leq i\leq 2\})\subset [1-\delta , 1]$$ and let $m$ be the length of the reduced word $W$. Let also $x_i = i/N, 0\leq i\leq N$ and $z = W(1/N)$.     
  
  \medskip
  
 By replacing the pair $(f,g)$ with $(f^{-1},g^{-1})$ if necessary, we may assume that $f(z)\geq z$. Then at least one of the following cases is valid:
 
  {\em Case 1}. $f(z)\leq gf(z)$;
  
  {\em Case 2}. $z\leq gf(z)$;
  
  {\em Case 3}. $gf(z)\leq z$.
  
  \medskip
  
  If Case 1 holds then we let $\alpha = f, \beta = g, z_0 = z$. If Case 1 does not hold but Case 2 holds, then we let $\alpha = gf, \beta = f, z_0 = z$. Finally, if Case 1 and Case 2 do not hold but Case 3 holds, then we let $\alpha = f^{-1}g^{-1}, \beta = g^{-1}, z_0 = gf(z)$. 
  
  \medskip
  
  In all the three cases, we will have $z_0\in [1-\delta , 1], z_0\leq z$, and $\alpha , \beta $ generate a free subsemigroup, and conditions of Lemma 1 are satisfied, i.e. we have $z_0\leq \alpha (z_0)\leq \beta \alpha (z_0)$. Moreover, we notice that $\mathrm{sup} _{0\leq x\leq 1}(|\alpha '(x)|+|\beta '(x)|)\leq M^2$, and the length of $W$ in the alphabet $\{\alpha , \beta , \alpha ^{-1}, \beta ^{-1}\}$ is at most $2m$. 
  
  \medskip
 
 Now, for every $n\in \mathbb{N}$, let $$\mathbb{S}_n = \{U(\alpha ,\beta )\beta \alpha W \ | \ U(\alpha ,\beta ) \ \mathrm{is \ a \ positive \ word \ in} \ \alpha , \beta  \ \mathrm{of \ length \ at \ most} \ n.\}$$ 
 
  \medskip

 Applying Lemma 1 to the pair $\{\alpha , \beta \}$ we obtain that $VW^{-1}(z_0) \geq z_0$ for all $V\in \mathbb{S}_n$.
 
 \medskip
   
  Then $|\mathbb{S}_n| \geq 2^n$. After achieving this inequality, we proceed as in the proof of Proposition 1 with just a slight change: there exists a sufficiently big $n$ such that the following two conditions are satisfied:
 
 \medskip
 
 (i) there exist $g_1, g_2\in \mathbb{S}_n$ such that $g_1\neq g_2$, and $$|g_1W(x_i)-g_2W(x_i)| < \frac{1}{\sqrt[2N]{2}^n} , 1\leq i\leq N-1.$$
 
 \medskip
 
 (ii) $M^{2m+4}(\theta _N)^n\frac{1}{\sqrt[2N]{2}^n} < \epsilon $.
 
 \bigskip
 
 Let $$h_1 = g_1W, h_2 = g_2W, y_i = W(x_i), z_i' = g_1(y_i), z_i'' = g_2(y_i), 1\leq i\leq N.$$ 
 
 \medskip
 
 Without loss of generality, we may also assume that $g_2(y_1) \geq g_1(y_1)$.
 
 \medskip
 
 Then for all $i\in \{1, \ldots , N-1\}$, we have 
 
 $$|h_1^{-1}h_2(x_i) - x_i| = |(g_1W)^{-1}(g_2W)(x_i) - x_i| =$$ \ $$|(g_1W)^{-1}(g_2W)(x_i) - (g_1W)^{-1}(g_1W)(x_i)| = |W^{-1}g_1^{-1}g_2(y_i) - W^{-1}g_1^{-1}g_1(y_i)|$$ \ $$ = |
 W^{-1}g_1^{-1}(z_i'') - W^{-1}g_1^{-1}(z_i')|$$
 
 \medskip
 
  Since $g_1, g_2\in \mathbb{S}_n$, by the Mean Value Theorem, we have $$|h_1^{-1}h_2(x_i) - x_i|\leq M^{2m+4}(\theta _N)^n|z_1'-z_1''| < M^{2m+4}(\theta _N)^n\frac{1}{\sqrt[2N]{2}^n}$$ By condition (ii), we obtain that $|h_1^{-1}h_2(x_i) - x_i| < \epsilon $. Then we have $|h_1^{-1}h_2(x) - x| < 2\epsilon $ for all $x\in [0,1]$. Since $\epsilon $ is arbitrary, we obtain that $\Gamma $ is not $C_0$-discrete. $\square $    
 
 \bigskip
 
 {\bf Proof of Theorem B.} Let $H$ be an arbitrary finitely generated subgroup of $[\Gamma , \Gamma ]$. If $H$ contains a non-abelian free subsemigroup then we are done by Theorem A. If $H$ does not contain a non-abelian free subsemigroup then by the result from {\bf [N2]} $H$ is virtually nilpotent. Then again by the result of {\bf [N2]}, $H$ is solvable of solvability degree at most $l(c)$. Since the natural number $l(c)$ depends only on $c$, and not on $H$, and since $H$ is an arbitrary finitely generated subgroup of $[\Gamma ,\Gamma ]$ we obtain that $[\Gamma , \Gamma ]$ is solvable of solvability degree at most $l(c)$. Hence $\Gamma $ is solvable with a solvability degree at most $l(c)+1$. $\square $ 
 
 \bigskip
 
 {\bf Proof of Theorem C.} Let again $H$ be an arbitrary finitely generated subgroup of $[\Gamma , \Gamma ]$. Again, if $H$ contains a non-abelian free subsemigroup then we are done by Theorem A. If $H$ does not contain a non-abelian free subsemigroup then by the result from {\bf [N2]} $H$ is virtually nilpotent. Then, by the result of Plante-Thurston ({\bf [PT]}), $H$ is virtually Abelian. Then, by the result of Szekeres ({\bf [S]}), $H$ is Abelian. Since $H$ is an arbitrary finitely generated subgroup of $[\Gamma ,\Gamma ]$, we conclude that $[\Gamma ,\Gamma ]$ is Abelian, hence $\Gamma $ is metaabelian. $\square $

  \vspace{3cm}
  
  {\bf R e f e r e n c e s:}

  [A] Akhmedov, A. \ On free discrete subgroups of Diff(I). \ {\em Algebraic and Geometric Topology}, \ {\bf vol.4}, (2010) 2409-2418. 
  
  [FF] Farb, B., Franks, J. \ Groups of homeomorphisms of one-manifolds III: Nilpotent subgroups. \ {\em Ergodic Theory and Dynamical Systems} \ {\bf 23} (2003), 1467-1484
  
 [N1] Navas, A. Sur les rapprochements par conjugasion en dimension 1 et classe $C^1$. \ {\em http://arxiv.org/pdf/1208.4815}.
  
 [N2] Navas, A. Growth of groups and diffeomorphisms of the interval. \ {\em Geom. Funct. Anal.} \ {\bf 18} 2008, {\bf no.3}, 988-1028. 
  
 [PT] Plante, J., Thurston, W. \ Polynomial growth in holonomy groups of foliations. \ {\em Comment. Math. Helv.} \ {\bf 51} (1976), 567-584. 
  
 [R] Raghunathan,M.S. {\em Discrete subgroups of semi-simple Lie groups}. \ Springer-Verlag, New York 1972. Ergebnisse der Mathematik und ihrer Grenzgebiete, Band 68. 
  
  [S] Szekeres, G. \ Regular iteration of real and complex functions. {\em Acta Math.} \ {\bf 100} (1958), 203-258

\end{document}